\documentclass[11pt]{article}
\usepackage{latexsym,amssymb,amsmath,amscd,amsthm,amsxtra}
\usepackage{mathrsfs}
\usepackage{epsfig}
\usepackage{pgf,tikz}
\usepackage{graphicx}
\usepackage{graphics}
\usepackage{setspace}
\usepackage[all,cmtip]{xy}
\newcommand{\counte}{section}

\newtheorem{prop}{\bf Proposition}[\counte]
\newtheorem{lemma}{\bf Lemma}[\counte]
\newtheorem{theorem}{\bf Theorem}[\counte]
\newtheorem{coro}{\bf Corollary}[\counte]
\newtheorem{example}{\bf Example}[\counte]

\newtheorem{conj}{\bf Conjecture}[\counte]
\newtheorem{remark}{\bf Remark}[\counte]

\def\Z{{\mathbb Z}}
\def\Q{{\mathbb Q}}

\author{Zhao Xu-an, zhaoxa@bnu.edu.cn\\ Gao Hongzhu, hzgao@bnu.edu.cn \\ School of Mathematical Sciences, Beijing Normal University\\Key Laboratory
of Mathematics and Complex Systems\\ Ministry of Education,
China, Beijing 100875 \\
Ruan Yangyang, ruanyy@amss.ac.cn\\ Institute of Mathematics, Chinese Academy of Sciences \\ China, Beijing 100190}

\title{The rational cohomology groups of the classifying spaces of Kac-Moody groups\thanks{The authors are supported by National Science Foundation of China, 12071034.}}

\date{}

\begin{document}

\maketitle
\begin{abstract}
In this paper, we compute the rational cohomology groups of the classifying space of a simply connected Kac-Moody group of infinite type. The fundamental principle is ``from finite to infinite''. That is, for a Kac-Moody group $G(A)$ of infinite type, the input data for computation are the rational cohomology of classifying spaces of parabolic subgroups of $G(A)$(which are of finite type), and the homomorphisms induced by inclusions of these subgroups. In some special cases, we can further determine the cohomology rings. %And we show that the algebraic difficulties which hinder us from the final computation results.
Our method also applies to study the mod $p$ cohomology of the classifying spaces of Kac-Moody groups.
\end{abstract}

\section{Introduction}
Let $A=(a_{ij})$ be an $n\times n$ integer matrix which satisfies:

(1) For each $i,a_{ii}=2$;

(2) For $i\not=j,a_{ij}\leq 0$;

(3) If $a_{ij}=0$, then $a_{ji}=0$.

\noindent then $A$ is called a Cartan matrix.

Let $h$ be the rational linear space with basis $\Pi^{\vee}=\{\alpha^{\vee}_1,\alpha^{\vee}_2,\cdots,\alpha^{\vee}_n\}$, and denote the dual basis of $\Pi^\vee$ in the dual space $h^*$ by $\{\omega_1,\omega_2,\cdots,\omega_n\}$, that is, $\omega_i(\alpha_j^{\vee})=\delta_{ij}$ for $1\leq i,j\leq n$. Let $\Pi=\{\alpha_1,\cdots,\alpha_n\}\subset h^*$ be given by $\alpha_j(\alpha^{\vee}_i)=a_{ij}$ for all $i,j$, then $\alpha_j=\sum\limits_{i=1}^n a_{ij}\omega_i$. $\Pi$ and $\Pi^\vee$ are called the simple root system and simple co-root system associated to the Cartan matrix $A$. $\alpha_i, \alpha_i^\vee, \omega_i,1\leq i\leq n$ are called respectively the simple roots, simple co-roots and fundamental dominant weights.

For each Cartan matrix $A$, there is a Lie algebra $g(A)$ associated to
$A$ called Kac-Moody Lie algebra, see Kac\cite{Kac_68}\cite{Kac_82} and Moody\cite{Moody_68}.

The Kac-Moody Lie algebra $g(A)$ is generated by elements $\alpha^\vee_i,e_i,f_i,1\leq i\leq n$ over $\mathbb{C}$,
with the defining relations:

(1) $[\alpha_i^{\vee},\alpha_j^{\vee}]=0, 1\leq i,j\leq n$;

(2) $[e_i,f_j]=\delta_{ij}\alpha_i^{\vee}, 1\leq i,j\leq n$;

(3) $[\alpha_i^{\vee},e_j]=a_{ij}e_j,[\alpha_i^{\vee},f_j]=-a_{ij}f_j, 1\leq i,j\leq n$;

(4) $\mathrm{ad}(e_i)^{-a_{ij}+1}(e_j)=0, 1\leq i\not=j\leq n$;

(5) $\mathrm{ad}(f_i)^{-a_{i j}+1}(f_j)=0, 1\leq i\not=j\leq n$.

Kac and Peterson\cite{Kac_Peterson_83}\cite{Kac_Peterson_84}\cite{Kac_85} constructed the Kac-Moody group $G(A)$ with Lie algebra $g(A)$.

A Cartan matrix $A$ is indecomposable if $A$ can not be written as a direct sum of two Cartan matrices
$A_1$ and $A_2$. Indecomposable Cartan matrices, their associated Kac-Moody Lie algebras, or Kac-Moody groups are divided into three types.

(1) Finite type, if $A$ is positive definite. In this case, $G(A)$ is just a simply connected complex semisimple Lie group.

(2) Affine type, if $A$ is positive semi-definite and has rank $n-1$.

(3) Indefinite type, otherwise.

The last two types of Cartan matrices or Kac-Moody groups are also called of infinite type. A decomposable Kac-Moody group is isomorphic to a product of Kac-Moody groups of these three types.

The Weyl group $W(A)$ associated to a Cartan matrix $A$ is the group generated by the Weyl reflections $\sigma_i:h^*\to h^*$ with respect to $\alpha_i^{\vee},1\leq i\leq n$, where $\sigma_i(\alpha)=\alpha-\alpha(\alpha_i^\vee)\alpha_i$. The action of $\sigma_i$ on the fundamental dominant weights is given by $\sigma_i(\omega_j)=\omega_j-\omega_j(\alpha_i^{\vee})\alpha_i=\omega_j-\delta_{ij}\alpha_i$. For details see Kac\cite{Kac_82}.

Let $S$ be the set $\{1,2,\cdots,n\}$. For each proper subset $I\subset S$, the matrix $A_I=(a_{ij})_{i,j\in I}$ is also a Cartan matrix. And there is a parabolic subalgebra $g_I(A)\subseteq g(A)$ generated by $e_i,i\in I$ and $\alpha_j^{\vee},f_j,1\leq j\leq n$. Corresponding to $g_I(A)$, there is a connected parabolic subgroup $G_I(A)$ of $G(A)$. For $I=S$, $G_I(A)=G(A)$. And for $I=\emptyset$, $G_I(A)=B(A)$, the Borel subgroup of $G(A)$.

Let $C(A)$ be the category whose objects are proper subsets $I$ of $S$ such that $A_I$ are of finite type, and morphisms are the inclusions of these proper subsets. There is a functor $F: C(A)\to Top$ which sends $I$ to the classifying space $BG_I(A)$ and sends $I\subset J$ to the map $BG_I(A)\to BG_J(A)$. By a theorem of Kichiloo\cite{Broto_Kitchloo_02}\cite{Kitchloo_17}, there is

\noindent {\bf Theorem. }
For a Cartan matrix $A$ of infinite type, the homotopy type of $BG(A)$ is the homotopy colimit of the functor $F$, i.e. $BG(A)\simeq \mathrm{hocolimit} F$.

This theorem enables us to compute the cohomology of $BG(A)$ from the cohomology of $BG_I(A)$ for those $A_I$ of finite type. It gives a principle that ``from finite to infinite''. By Bousfield and Kan\cite{Bousfield_Kan_72}, there is a multiplicative spectral sequence with $E_2$-term given by the derived limits of the functor $H^*\circ F$ which converges to the cohomology of $BG(A)$. This spectral sequence is used in Kitchloo\cite{Kitchloo_17} to study the cohomology of the classifying spaces of Kac-Moody groups. But the $E_2$-term of Bousfield-Kan spectral sequence is hard to compute. So in this paper we introduce a homotopy simplification procedure. In fact we show that there are spaces $X; X_{I_1},X_{I_2},\cdots,X_{I_d}$ for certain $d\geq 2$ such that $X=X_{I_1}\cup X_{I_2}\cup \cdots \cup X_{I_d}$, $X$ has the same homotopy type as $BG(A)$, and for $1\leq i\leq d$, $X_{I_i}$ has the same homotopy type as $BG_{I_i}(A)$. These spaces $X$ and $X_{I_i}, 1\leq i\leq d$ are constructed explicitly in section 2. As a consequence we can use Mayer-Vietoris spectral sequence of $X=X_{I_1}\cup X_{I_2}\cup \cdots \cup X_{I_d}$ to compute the cohomology groups of $BG(A)$.

The contents of this paper are as follows. In section 2, we simplify the homotopy colimit of the classifying space $BG(A)$ and construct the spaces $X; X_{I_1}, X_{I_2},\cdots, X_{I_d}$. In section 3 we use the Mayer-Vietoris spectral sequence to compute the rational cohomology groups of $BG(A)$. We will explain how to compute the $E_2$-term of this spectral sequence. And we prove that for rational coefficients the spectral sequence collapses at $E_2$-term. So it gives the cohomology groups of $BG(A)$. In section 4 we give an algorithm to compute the cohomology groups and list the computation results for the cases of $n=3$ and $4$. In section 5 we give some results about the rational polynomial invariants of Weyl groups $W(A)$, which appear naturally in the M-V spectral sequence. In section 6 we discuss some applications of our results and give some conjectures on the cohomology of classifying spaces of Kac-Moody groups.

\section{The simplification of the homotopy colimit}
Let us recall the definition of homotopy colimit. See Dugger\cite{Dugger_08} for details.

Let $C$ be a small category and $F: C\rightarrow Top$ be a functor from $C$ to the category of topological spaces. Then we have a simplicial space $\mathrm{srep}(F)$ which is called the simplicial replacement of $F$
with $\mathrm{srep}(F)_n =\coprod\limits_{i_0\rightarrow i_1\rightarrow\cdots \rightarrow i_n} F(i_0)$, the coproduct ranges over chains of composable morphisms in $C$. The face and degeneracy maps are defined as follows. If $\sigma: i_0\rightarrow i_1\rightarrow\cdots \rightarrow i_n$ is a chain, then by deleting $i_j, j\geq0$ we obtain a chain $\sigma(j)$ of $n-1$ composable maps. When $j >0$, the map $d_j: \mathrm{srep}(F)_n \rightarrow \mathrm{srep}(F)_{n-1}$ sends the summand $F(i_0)$ corresponding to the chain $\sigma$ to the identical copy of $F(i_0)$ in $\mathrm{srep}(F)_{n-1}$ indexed by $\sigma(j)$. When $j = 0, d_0: \mathrm{srep}(F)_n \rightarrow  \mathrm{srep}(F)_{n-1}$ sends the summand $F(i_0)$ corresponding to the chain $\sigma$ to the summand $F(i_1)$ corresponding to $\sigma(0)$ and
the map we use here is $F(i_0)\rightarrow F(i_1)$ given by the first map $i_0\to i_1$ in $\sigma$. For the degeneracy maps each $s_j$ sends the summand $F(i_0)$ corresponding to the chain $\sigma$ to the identical summand $F(i_0)$ corresponding to the chain in which one has inserted the identity map $i_j \rightarrow i_j$. The homotopy colimit of the functor $F: C\rightarrow Top$ is defined as the geometric
realization of its simplicial replacement $\mathrm{srep}(F)$.

For a simplicial complex $L$ with vertex set $S=\{1,2,\cdots,n\}$, there is a small category $C(L)$ with objects the simplices of $L$. And the morphisms of $C(L)$ are given by inclusions of simplices. We called the category $C(L)$ the simplicial category associated to $L$. In a simplicial category $C(L)$ all the automorphisms are identities, and for two objects there is at most one morphism between them. Given a Cartan matrix $A$, the category $C(A)$ defined in the introduction is a simplicial category since if $A_J$ is of finite type, then for any $I\subset J$, $A_I$ is also of finite type. And for each $i\in S$, $\{i\}\in C(A)$ because the rank 1 Cartan matrix is of finite type. Therefore for a Kac-Moody group $G(A)$, it is convenient to regard the category $C(A)$ as a proper simplicial sub-complex of the $n-1$-simplex $\Delta_{n-1}$ with vertex set $S=\{1,2,\cdots,n\}$, and regard an object $I\in C(A)$ as a simplex.

Let $C$ be the simplicial category associated to a simplicial complex $L$ and $F: C\to Top$ be a functor, For an object $I\in C$, let $C_I$ be the full subcategory of $C$ whose objects are $I$ and faces of $I$. We define $X_I=\mathrm{hocolimit} F|_{C_I}$, the homotopy colimit of $F$ restricted to $C_I$. Since $C_I$ has a terminal object $I$ we know $X_I\simeq F(I)$.
The inclusion $I\subset J$ induces the map $i_{IJ}: X_I\rightarrow X_J$. As a summary, we have

\begin{lemma}

1. For $I\subset J$, the map $i_{IJ}: X_I\rightarrow X_J$ is an inclusion(In fact a cofibration).

2. For each $I\in C$, $X_I$ is a simplicial subspace of $X=\mathrm{hocolimit} F$.

3. $X_I\cap X_J=X_{I \cap J}$ for $I,J\in C$.

\end{lemma}

This lemma can be checked by the definition of homotopy colimits. It can be regarded as a homotopy replacement result for the functor $F$. By applying it to the functor $F: C(A)\rightarrow Top$, we can replace the map $BG_I(A)\rightarrow BG_J(A)$ by a cofibration $X_I\subset X_J$, and the classifying space $BG_I(A)$ by $X_I$. Then the homotopy colimit of $F$ is replaced by a union of subspaces $X_I$ in $X$.

%We call a simplex $\sigma$ of a simplicial complex $L$ which is not a face of any other simplex, a maximal simplex of $L$.
We call a simplex $\sigma$ of a simplicial complex $L$ a maximal simplex if it is not a face of any other simplex. It is obvious that the geometric realizations of a simplicial complex is the union of geometric realization of its maximal simplices. For a finite simplicial category $C=C(L)$, the homotopy colimit of a functor $F:C\rightarrow Top$ has a simple description.

\begin{theorem}
Let $L$ be a simplicial complex contained in $\Delta_{n-1}$ and $I_1,I_2,\cdots I_d$ be the set of maximal simplices of $L$, then for a functor $F: C(L)\rightarrow Top$, we have $X=\mathrm{hocolimit} F = \bigcup\limits_{1\leq i\leq d} X_{I_i}$.
\end{theorem}

Combining Theorem 2.1 with Kitchloo's theorem, we have

\begin{coro}
Let $A$ be a Cartan matrix of infinite type, $I_1,I_2,\cdots,I_d\in C(A)$ be the maximal objects of $C(A)$, then $BG(A)\simeq \bigcup\limits_{1\leq i\leq d} X_{I_i}$.
\end{coro}

\begin{example}
For a Cartan matrix $A$ of rank $3$, up to permutation of $1,2,3$, the simplicial category $C(A)$ corresponding to $A$ is one of the following $4$ types.

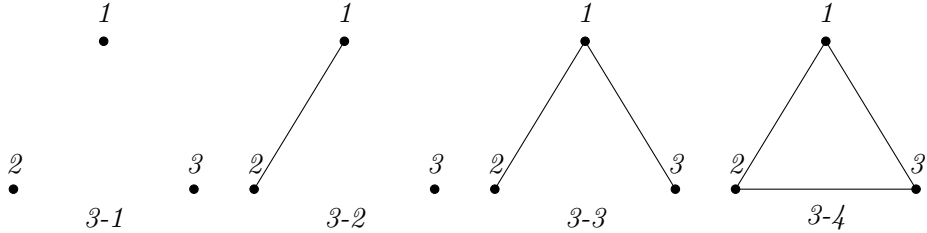
\begin{center}
\begin{tikzpicture}[scale=0.8]

{\filldraw [black] (18-12,-0.73) circle (2pt);}
{\filldraw [black] (21-12,-0.73) circle (2pt);}
{\filldraw [black] (19.5-12,1.73) circle (2pt);}

\draw (19.5-12,-1.2)  node{3-1};

{\filldraw [black] (18-8,-0.73) circle (2pt);}
{\filldraw [black] (21-8,-0.73) circle (2pt);}
{\filldraw [black] (19.5-8,1.73) circle (2pt);}

{\draw (18-8,-0.73) -- +(1.5, 2.46);}
\draw (19.5-8,-1.2)  node{3-2};

{\filldraw [black] (18-4,-0.73) circle (2pt);}
{\filldraw [black] (21-4,-0.73) circle (2pt);}
{\filldraw [black] (19.5-4,1.73) circle (2pt);}

{\draw (18-4,-0.73) -- +(1.5, 2.46);}
{\draw (21-4,-0.73) -- +(-1.5, 2.46);}
\draw (19.5-4,-1.2)  node{3-3};

{\filldraw [black] (18,-0.73) circle (2pt);}
{\filldraw [black] (21,-0.73) circle (2pt);}
{\filldraw [black] (19.5,1.73) circle (2pt);}

{\draw (18,-0.73) -- +(1.5, 2.46);}
{\draw (21,-0.73) -- +(-1.5, 2.46);}
{\draw (18,-0.73) -- +(3, 0);}
\draw (19.5,-1.2)  node{3-4};

\draw [black] (18-12,-0.3)  node{2};
\draw [black] (18-8,-0.3)  node{2};
\draw [black] (18-4,-0.3)  node{2};
\draw [black] (18,-0.3)  node{2};

\draw [black] (21-12,-0.3)  node{3};
\draw [black] (21-8,-0.3)  node{3};
\draw [black] (21-4,-0.3)  node{3};
\draw [black] (21,-0.3)  node{3};

\draw [black] (19.5-12,2.2) node{1};
\draw [black] (19.5-8,2.2) node{1};
\draw [black] (19.5-4,2.2) node{1};
\draw [black] (19.5,2.2) node{1};

\end{tikzpicture}
\end{center}

In the four cases the maximal objects of $C(A)$ and the homotopy type of $BG(A)$ are

(1). $\{1\},\{2\},\{3\}, \ \ BG(A)\simeq (X_1\cup_{X_{\emptyset}} X_2) \cup_{X_{\emptyset}} X_3$;

(2). $\{1,2\}, \{3\}, \ \ BG(A)\simeq X_{12}\cup_{X_{\emptyset}} X_3$;

(3). $\{1,2\}, \{1,3\}, \ \ BG(A)\simeq X_{12}\cup_{X_1} X_{13}$;

(4). $\{1,2\},\{1,3\}, \{2,3\}, \ \ BG(A)\simeq (X_{12}\cup_{X_1} X_{13})\cup_{X_{2}\cup_{X_{\emptyset}} X_3} X_{23} $.
\end{example}

%\begin{prop}
%For each sub-complex $C$ of $\Delta_{n-1}$ which contain all the vertices set $S$, there is a Cartan matrix $A$ such that $C(A)=C$.
%\end{prop}

%{\bf Proof: }Let $C$ be a sub-complex of $\Delta_{n-1}$ which contain all the vertices set. Let the maximal of $C$ be $I_1,\cdots,I_n$, then define a Cartan matrix $A$ with $a_{ij}=0$ if $i,j\in I_t$ for $1\leq t\leq d$. for  $S$$simplical
%We do not know whether each admissible category $C$ is obtained from a Kac-Moody group $K(A)$.
\section{Mayer-Vietoris spectral sequence}

In section 2 we write the homotopy type of the classifying space $BG(A)$ as a union of subspaces $X_{I_1},X_{I_2},\cdots,X_{I_d}$ which are homotopy equivalent to $BG_{I_1}(A), BG_{I_2}(A), \cdots, BG_{I_d}(A)$ respectively. Now we can use the Mayer-Vietoris exact sequences to compute the cohomology groups of $BG(A)$ inductively. The computation procedure can be organized into an M-V spectral sequence. The spectral sequence $E_r,r\geq 0$ is constructed from the bi-complex $E_0^{r,s}=\bigoplus\limits_{i_1<i_2<\cdots <i_{r+1}} S^s_{\mathcal{U}}(X_{I_{i_1}}\cap X_{I_{i_2}}\cap\cdots\cap X_{I_{i_{r+1}}})$ with horizontal differential $\delta$ and vertical differential $d$. The $E_1$-term and $E_2$-term of this spectral sequence are given by $H^d(E_0)$ and $H^\delta E_1=H^\delta H^d (E_0)$. This spectral sequence converges to $H^*(X)$. See Bott and Tu\cite{Bott_Tu_82} for details. By applying it to $X_{I_1}\cup X_{I_2}\cup \cdots\cup X_{I_d}$, we have

\begin{theorem}
For arbitrary abelian group $G$, there is an M-V spectral sequence with $E_1^{r,s}=\bigoplus\limits_{i_1<i_2<\cdots <i_{r+1}} H^s(X_{I_{i_1}}\cap X_{I_{i_2}}\cap\cdots\cap X_{I_{i_{r+1}}};G)$, which converges to the cohomology $H^*(BG(A);G)$.
\end{theorem}

As a consequence, with coefficient group a field $k$, we have $H^*(BG(A);k)\cong E_\infty^{0,*}\oplus E_\infty ^{1,*-1}\oplus \cdots \oplus E^{*,0}_\infty$. From the expression of $E_1^{r,s}$ it is obvious that for $r\geq d$, $E_1^{r,s}=0$. So the spectral sequence is concentrated on the range of $0\leq r\leq d-1$.
Let $E^{r}_{\infty}=\bigoplus\limits_{s=0}^\infty E_\infty^{r,s}$ for $r\geq 0$, then $E^r_\infty =0$ for $r\geq d$ and we have $H^*(BG(A);k)\cong \bigoplus\limits_{r=0}^{d-1} \Sigma^r E^r_{\infty}$ as a graded $k$-vector space, where $\Sigma^r$ lifts the degree of a graded vector space by $r$.

In the following all cohomology groups are of rational coefficients.

The input data to compute the rational cohomology $H^*(BG(A))$ come from the functor $H^*\circ F: C(A)\to \mathbb{Q}Alg$ which sends $I\in C(A)$ to $H^*(BG_I(A))$ and sends a morphism $I\subset J$ to $H^*(BG_J(A))\to H^*(BG_I(A))$. The ring $H^*(BG_I(A))$ and homomorphism $H^*(BG_J(A))\rightarrow H^*(BG_I(A))$ can be computed by the following theorem of Borel, see \cite{Borel_53_1}.

\begin{theorem}
Let $K$ be a connected compact Lie group with Weyl group $W(K)$, and $T$ be a maximal torus subgroup of $K$, then the cohomology homomorphism $Bi^*:H^*(BK;\Q)\rightarrow H^*(BT;\Q)$ induced by the inclusion $T\to K$ is injective. And its image is exactly the invariants $H^*(BT;\Q)^{W(K)}$ of Weyl group.
\end{theorem}

For a parabolic subgroup $G_I(A)$ of finite type, the unitary form of $G_I(A)$ is a compact Lie subgroup $K_I(A)$ of the unitary form $K(A)$ of $G(A)$. And the inclusion $K_I(A)\subset G_I(A)$ is a homotopy equivalence. Hence for the computation of the cohomology we can replace the group $G_I(A)$ by it's unitary form $K_I(A)$. Since the Weyl group of $G_I(A)$(or $K_I(A)$) is $W_I(A)$, which is generated by Weyl reflection $\sigma_{i},i\in I$, the ring $H^*(BG_I(A))\cong H^*(BK_I(A))$ can be identified with $H^*(BT)^{W_I(A)}$ by Borel's theorem. For convenience we write $P$ for the polynomial algebra $H^*(BT)$, which is naturally identified with $\Q[w_1,\cdots,w_n]$, where $w_i,1\leq i\leq n$ are the fundamental dominant weights. And for each $I\subset S$ we write $P_I$ for the polynomial subalgebra $H^*(BT)^{W_I(A)}$. If we identify $H^*(BG_I(A))$ with $P_I$, then the homomorphism $Bi_{IJ}^*:H^*(BG_J(A))\rightarrow H^*(BG_I(A))$ for $I\subset J$ can be identified with the inclusion $P_J=H^*(BT)^{W_J(A)}\subset H^*(BT)^{W_I(A)}=P_I$. An easy but useful fact is that $P_{I_1}\cap P_{I_2}=P_{I_1\cup I_2}$.

%We denote the homotopy replacement of $BK_I(A)$ by $X_I$ as last section. So $X_{\emptyset}\simeq BT$. Let $P$ be the ring $H^*(BT)$ In the following for simplicity for each $I\in \{1,2,\cdots,n\}$ we denote %$H^*(BT)^{W_I(A)}$ by $P_I$ and $P_{\emptyset}=H^*(BT)$.

%Let $P$ be the ring $H^*(BT)$, for each $I\in S$, we have a sub-ring $P_I$ which is the polynomial invariants of Weyl group $W_I(A)\subset W(A)$. And we set $P_{\emptyset}=P$.

Let $E^r_1=\bigoplus\limits_{s=0}^\infty E_1^{r,s}$ for $0\leq r\leq d-1$. From the above discussion we have $$E_1^{r}\cong \bigoplus\limits_{s=0}^\infty \bigoplus\limits_{i_1<i_2<\cdots <i_{r+1}} H^s(X_{I_{i_1}}\cap X_{I_{i_2}}\cap\cdots\cap X_{I_{i_{r+1}}};\Q)\cong \bigoplus\limits_{i_1<i_2<\cdots <i_{r+1}} P_{I_{i_1}\cup I_{i_2}\cup \cdots \cup I_{ i_{r+1}}}.$$ The $E_2$-term of the M-V spectral sequence can be computed from $E_1$-term as the homology of the following chain complex $0\rightarrow  E_{1}^0\stackrel{\delta_0}\rightarrow E_{1}^1\stackrel{\delta_1}\rightarrow E_{1}^2 \stackrel{\delta_2}\rightarrow \cdots \stackrel{\delta_{d-2}}\rightarrow E_{1}^{d-1} \rightarrow 0$.

The following is an example.
%for the case $X=X_{I_1}\cup X_{I_2}\cup X_{I_3}\cup X_{I_4}$ with $d=4$.

\begin{example}

For the case $X=X_{I_1}\cup X_{I_2}\cup X_{I_3}\cup X_{I_4}$ with $d=4$, the chain complex can be depicted as follows.
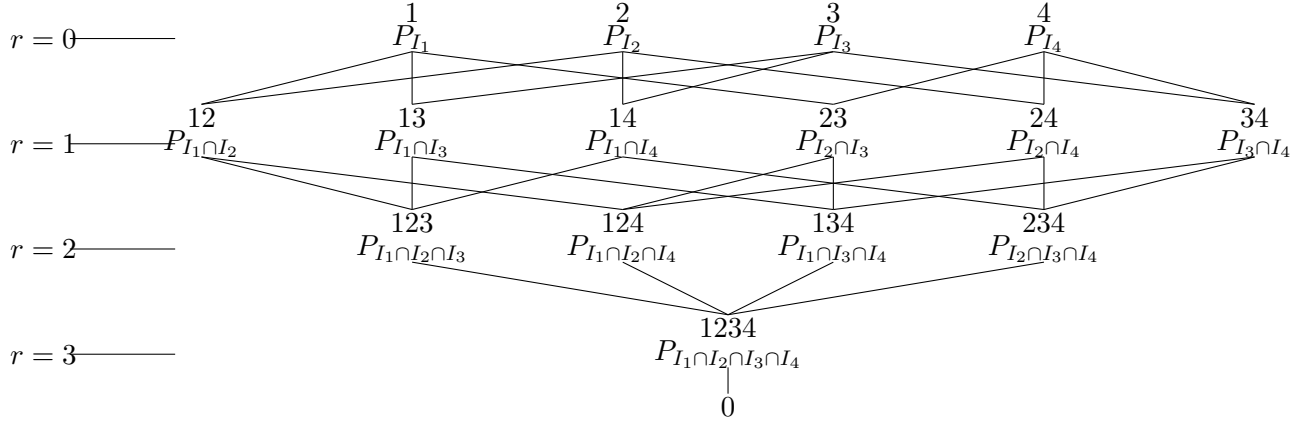
\begin{center}
\begin{tikzpicture}[scale=0.7]

\draw (-7,0)  node{$r=0$};  \draw (-8+1.5,0)  -- +(2,0) ;
\draw (-7,-2)  node{$r=1$}; \draw (-8+1.5,-2) -- +(2,0)  ;
\draw (-7,-4)  node{$r=2$}; \draw (-8+1.5,-4) -- +(2,0)  ;
\draw (-7,-6)  node{$r=3$}; \draw (-8+1.5,-6) -- +(2,0)  ;

\draw (0,0)  node{$P_{I_1}$};
\draw (4,0)  node{$P_{I_2}$};
\draw (8,0)  node{$P_{I_3}$};
\draw (12,0)  node{$P_{I_4}$};
\draw (0,0.25*2)  node{$1$};
\draw (4,0.25*2)  node{$2$};
\draw (8,0.25*2)  node{$3$};
\draw (12,0.25*2)  node{$4$};

\draw (0,-0.25)  -- (-4,-2+0.75);
\draw (0,-0.25)  -- (0,-2+0.75);
\draw (0,-0.25)  -- (8,-2+0.75);

\draw (4,-0.25) --  (-4,-2+0.75) ;
\draw (4,-0.25) --  (4,-2+0.75) ;
\draw (4,-0.25) --  (12,-2+0.75) ;

\draw (8,-0.25) --  (0,-2+0.75) ;
\draw (8,-0.25) --  (4,-2+0.75) ;
\draw (8,-0.25) --  (16,-2+0.75) ;

\draw (12,-0.25) --  (8,-2+0.75) ;
\draw (12,-0.25) --  (12,-2+0.75) ;
\draw (12,-0.25) --  (16,-2+0.75) ;

\draw (-4,-2)  node{$P_{I_1\cap I_2}$};
\draw (0,-2)  node{$P_{I_1\cap I_3}$};
\draw (4,-2)  node{$P_{I_1\cap I_4}$};
\draw (8,-2)  node{$P_{I_2\cap I_3}$};
\draw (12,-2)  node{$P_{I_2\cap I_4}$};
\draw (16,-2)  node{$P_{I_3\cap I_4}$};
\draw (-4,-2+0.25*2)  node{$12$};
\draw (0,-2+0.25*2)  node{$13$};
\draw (4,-2+0.25*2)  node{$14$};
\draw (8,-2+0.25*2)  node{$23$};
\draw (12,-2+0.25*2)  node{$24$};
\draw (16,-2+0.25*2)  node{$34$};

\draw (-4,-2-0.25)  --  (0,-4+0.75);
\draw (-4,-2-0.25)  --  (4,-4+0.75);

\draw (0,-2-0.25)   --  (0,-4+0.75);
\draw (0,-2-0.25)   --  (8,-4+0.75);

\draw (4,-2-0.25)   --  (0,-4+0.75);
\draw (4,-2-0.25)   --  (12,-4+0.75);

\draw (8,-2-0.25)   --  (4,-4+0.75);
\draw (8,-2-0.25)   --  (8,-4+0.75);

\draw (12,-2-0.25)  --  (4,-4+0.75);
\draw (12,-2-0.25)  --  (12,-4+0.75);

\draw (16,-2-0.25)  --  (8,-4+0.75);
\draw (16,-2-0.25)  --  (12,-4+0.75);

\draw (0,-4)  node{$P_{I_1\cap I_2\cap I_3}$};
\draw (4,-4)  node{$P_{I_1\cap I_2\cap I_4}$};
\draw (8,-4)  node{$P_{I_1\cap I_3\cap I_4}$};
\draw (12,-4)  node{$P_{I_2\cap I_3\cap I_4}$};
\draw (0,-4+0.5)  node{$123$};
\draw (4,-4+0.5)  node{$124$};
\draw (8,-4+0.5)  node{$134$};
\draw (12,-4+0.5)  node{$234$};

\draw (0,-4-0.25) -- (6,-6+0.75);
\draw (4,-4-0.25) -- (6,-6+0.75);
\draw (8,-4-0.25) -- (6,-6+0.75);
\draw (12,-4-0.25) -- (6,-6+0.75);

\draw (6,-6) node{$P_{I_1\cap I_2\cap I_3\cap I_4}$};
\draw (6,-6+0.5)  node{$1234$};

\draw (6,-6-0.25) -- ((6,-7+0.25);

\draw (6,-7)  node{$0$};

\end{tikzpicture}
\end{center}

\end{example}

In this diagram we omit the sign $\oplus$ between the adjacent items in each row. For the description of the differential see Bott and Tu\cite{Bott_Tu_82} or Example 4.1 in the next section.

For $I\in C(A)$, by Borel's theorem $P_I$ is concentrated on the even degree. So we have $E^{r,s}_1=\{0\}$ for $s$ odd. Then we have $E_k^{r,s}=\{0\}$ for $s$ odd, $k\geq 1$. We need this observation to prove the following theorem.

%\begin{lemma}The $E_2$-term of the M-V spectral sequence satisfies

%1. $E^{r,s}_2=\{0\}$ for $r\geq d$.

%2. $E^{r,s}_2=\{0\}$ for $s$ odd.

%\end{lemma}

\begin{theorem}
The rational M-V spectral sequence for $BG(A)$ collapses at $E_2$-term.
\end{theorem}

\noindent {\bf Proof: }By Kitchloo\cite{Kitchloo_17}, let $l$ be a prime number so that Weyl group $W(A)$ contains no elements of order $l$, then there is an unstable Adams map $\psi$ on $BG(A)$, and it induces the Adams map $\psi_I$ on $BG_I(A)$ for all $I\in C(A)$. The cohomology endomorphism $\psi_I^*$ acting on $H^{2q}(BG_I(A))$ for $G_I(A)$ of finite type is the multiplication by $l^q$. The map $\psi$ induces an endomorphism on the M-V spectral sequence. For even integer $2k$, $d^{r,s}_{2k}=0$ since the domain or the range of $d^{p,q}_{2r}$ is trivial by the observation $E_{2k}^{r,s}=\{0\}$ for $s$ odd, $k\geq 1$. For odd integer $2k-1$, if $x\in E_{2k-1}^{r,2s}$, since $\psi^*$ is commutative with the differential of the M-V spectral sequence, we have $d_{2k-1}\psi^*(x)=\psi^* d_{2k-1}(x)$, hence we get $l^{s}d_{2k-1}(x)=l^{s-k+1}d_{2k-1}(x)$. As a result if $k>1$, we have $d_{2k-1}x=0$. This shows that we always have $d_k(x)=0$, for $k\geq 2$. Therefore the M-V spectral sequence collapses.

\begin{remark}
There are similar M-V spectral sequences for all coefficient groups other than $\Q$ which converge to the associated cohomology group $H^*(BG(A))$. But in general they do not collapse at $E_2$-term.
\end{remark}

The explicit computation of the $E_2$-term of the rational M-V spectral sequence will be discussed in the next section.

\section{The computation of rational cohomology groups of $BG(A)$}
In this section we compute the rational cohomology groups of the classifying space $BG(A)$. Since for rational coefficient the corresponding M-V spectral sequence collapses at $E_2$-term, the $E_2$-term of the spectral sequence gives the rational cohomology groups of $BG(A)$.

The following is an example of explicit computation.

\begin{example}
We demonstrate the computation procedure for the M-V spectral sequence of $X=X_{I_1}\cup X_{I_2}\cup X_{I_3}\cup X_{I_4}$. The cochain complex $E_1: 0\rightarrow  E_{1}^0\stackrel{\delta_0}\rightarrow E_{1}^1\stackrel{\delta_1}\rightarrow E_{1}^2 \stackrel{\delta_2}\rightarrow \cdots \stackrel{\delta_{d-2}}\rightarrow E_{1}^{d-1} \rightarrow 0$ is illustrated as in Example 3.1. Here we compute the cohomology $E_2=H^*(E_1)=\bigoplus\limits_{r= 0}^{d-1} H^r(E_1)$.

{\bf Step 1:} It is easy to see that $H^0(E_1)\cong P_{I_1}\cap P_{I_2}\cap P_{I_3}\cap P_{I_4}=P_{I_1\cup I_2\cup I_3\cup I_4}$.

{\bf Step 2:} For $H^1(E_1)$, the element in $\ker \delta_1$ has the form $(v_{12},v_{13},v_{14},v_{23},v_{24},v_{34})$ with $v_{ij}\in P_{I_i\cap I_j},1\leq i<j\leq 4$ and satisfies
\begin{equation}v_{12}-v_{13}+v_{23}=0.\end{equation}
\begin{equation}v_{12}-v_{14}+v_{24}=0.\end{equation}
\begin{equation}v_{13}-v_{14}+v_{34}=0.\end{equation}
\begin{equation}v_{23}-v_{24}+v_{34}=0.\end{equation}

The image of $\delta_0$ is generated by the elements of the forms $\delta_0(v_i)$ with $v_i\in P_{I_i},1\leq i\leq 4$, which is represented by the row vectors of the following matrix.

\begin{center}

$\left(
  \begin{array}{cccccc}
    12 & 13 & 14 & 23 & 24 & 34 \\
    v_1 & v_1 & v_1 & 0 & 0 & 0 \\
    -v_2 & 0 & 0 & -v_2 & -v_2 & 0 \\
    0 & -v_3 & 0 & -v_3 & 0 & v_3 \\
    0 & 0 & 0 & -v_4 & -v_4 & -v_4 \\
  \end{array}
\right)$
\end{center}

By the equations $v_{12}-v_{13}+v_{23}=0$ and $v_{12}-v_{14}+v_{24}=0$, we have $v_{12}\in P_{I_1\cap I_2}\cap (P_{I_1\cap I_3}+P_{I_2\cap I_3})\cap (P_{I_1\cap I_4}+P_{I_2\cap I_4})$. We compute $H^1(E_1)$ by constructing an epimorphism
$$\phi_1: H_1(E_1)\to P_{I_1\cap I_2}\cap (P_{I_1\cap I_3}+P_{I_2\cap I_3})\cap (P_{I_1\cap I_4}+P_{I_2\cap I_4})-(P_{I_1}+P_{I_2})$$ $$[(v_{12},v_{13},v_{14},v_{23},v_{24},v_{34})]\to v_{12}+P_{I_1}+P_{I_2}.$$
Here we use $A-B$ to denote the quotient of $A$ with respect to subspace $B$ for convenience. It is easy to see that $\phi_1$ is well defined. And as a rational vector space $H_1(E_1)\cong \mathrm{im}\phi_1\oplus \ker \phi_1$.

Now we compute $\ker\phi_1$. For $[(v_{12},v_{13},v_{14},v_{23},v_{24},v_{34})]\in \ker\phi_1$, we must have $v_{12}\in P_{I_1}+P_{I_2}$. Therefore there exist $v_{1}\in P_{I_1},v_{2}\in P_{I_2}$ such that $v_{12}=v_{1}-v_{2}$. Substituting $v_{12}$ with $v_1-v_2$ in Equations 1 and 2, we get $v_{13}-v_{1}=v_{23}-v_{2}$ and $v_{14}-v_{1}=v_{24}-v_{2}$. This means that $v_{13}-v_{1}\in P_{I_1\cap I_3}\cap P_{I_2\cap I_3}=P_{(I_1\cup I_2)\cap I_3}$ and $v_{14}-v_{1}\in P_{I_1\cap I_4}\cap P_{I_2\cap I_4}=P_{(I_1\cup I_2)\cap I_4}$. By Equation 3, $v_{13}-v_{14}+v_{34}=0$, we get $(v_{13}-v_1)-(v_{14}-v_1)+v_{34}=0$. So $v_{13}-v_1\in P_{(I_1\cup I_2)\cap I_4}+P_{I_3\cap I_4}$. Therefore $v_{13}-v_1\in P_{(I_1\cup I_2)\cap I_3} \cap (P_{(I_1\cup I_2)\cap I_4}+P_{I_3\cap I_4})$.

We define a homomorphism $$\phi_2: \ker \phi_1\rightarrow P_{(I_1\cup I_2)\cap I_3} \cap (P_{(I_1\cup I_2)\cap I_4}+P_{I_3\cap I_4})-(P_{I_1\cup I_2}+P_{I_3})$$
$$[(v_{12},v_{13},v_{14},v_{23},v_{24},v_{34})]\to v_{13}-v_1+P_{I_1\cup I_2}+P_{I_3}.$$

It is easy to check that this homomorphism is well defined. Note that we need the $P_{I_1\cup I_2}$-item to eliminate the effect of indeterminacy of $v_1$, and the $P_{I_3}$-item to send $[(0,-v_3,0,-v_3,0,v_3)]\in \ker \delta_0$ to $0$. By definition $\phi_2$ is an epimorphism and we have $\ker \phi_1\cong \mathrm{im}\phi_2\oplus \ker \phi_2.$

For $[(v_{12},v_{13},v_{14},v_{23},v_{24},v_{34})]\in \ker\phi_2$, we already have $v_1\in P_{I_1}$ and $v_2\in P_{I_2}$ such that $v_{12}=v_1-v_2$ for $\ker \phi_2\subset \ker \phi_1$, and there are $\alpha\in P_{I_1\cup I_2}$ and $v_3\in P_{I_3}$ such that $v_{13}-v_1=\alpha-v_3$. So we have $v_{12}=(v_1+\alpha)-(v_2+\alpha)$ and $v_{13}=(v_1+\alpha)-v_3$. This means that if we choose $v'_1=v_1+\alpha$ and $v'_2=v_2+\alpha$ instead of $v_1,v_2$ in the previous step, then we will have $v_{12}=v'_1-v'_2$ and $v_{13}=v'_1-v_3$. For simplicity we still use $v_1,v_2$ to denote the newly chosen elements $v'_1,v'_2$, then $\alpha$ disappears and we have $v_{12}=v_1-v_2$ and $v_{13}=v_1-v_3$. Combining them with Equation 1, we can solve $v_{23}=v_2-v_3$. By Equations 2, 3 and 4, we get $v_{14}-v_1=v_{24}-v_2=v_{34}-v_3$. These equations imply that $v_{14}-v_1\in P_{I_1\cap I_4}\cap P_{I_2\cap I_4}\cap P_{I_3\cap I_4}=P_{(I_1\cup I_2\cup I_3)\cap I_4}$.

Now we define an epimorphism $$\phi_3: \ker \phi_2\rightarrow P_{(I_1\cup I_2\cup I_3)\cap I_4}-(P_{I_1\cup I_2\cup I_3}+P_{I_4})$$
$$[(v_{12},v_{13},v_{14},v_{23},v_{24},v_{34})]\to v_{14}-v_1+P_{I_1\cup I_2\cup I_3}+P_{I_4}.$$

For $[(v_{12},v_{13},v_{14},v_{23},v_{24},v_{34})]\in \ker \phi_3$, we have $v_{12}=v_{1}-v_2,v_{13}=v_1-v_3, v_{23}=v_2-v_3$ and $v_{14}-v_1=\beta-v_4$ for certain $v_i\in P_{I_i},1\leq i\leq 4$ and $\beta\in P_{I_1\cup I_2\cup I_3}$. If we replace $v_1,v_2,v_3$ in the previous step by $v'_1=v_1+\beta,v'_2= v_2+\beta,v'_3=v_3+\beta$, then $\beta$ disappears, and we have $v_{14}-v_1=-v_4$. As a result we get $v_{14}=v_1-v_4, v_{24}=v_2-v_4$ and $v_{34}=v_3-v_4$ from the Equations 2, 3 and 4. Hence $(v_{12},v_{13},v_{14},v_{23},v_{24},v_{34})$ is in $\mathrm{im} \delta_0$ and $\ker \phi_3=0$. Hence $\ker \phi_2\cong \mathrm{im} \phi_3\oplus \ker \phi_3\cong \mathrm{im} \phi_3$.

So we have

$H^1(E_1)\cong \mathrm{im} \phi_1\oplus \mathrm{im} \phi_2\oplus \mathrm{im} \phi_3 $

$\cong P_{I_1\cap I_2}\cap (P_{I_1\cap I_3}+P_{I_2\cap I_3})\cap (P_{I_1\cap I_4}+P_{I_2\cap I_4})-(P_{I_1}+P_{I_2})$

$ \oplus P_{(I_1\cup I_2)\cap I_3} \cap (P_{(I_1\cup I_2)\cap I_4}+P_{I_3\cap I_4})-(P_{I_1\cup I_2}+P_{I_3})\oplus P_{(I_1\cup I_2\cup I_3)\cap I_4}-(P_{I_1\cup I_2\cup I_3}+P_{I_4})$

{\bf Step 3:} For $H^2(E_1)$, the element in $\ker \delta_2$ has the form $(v_{123},,v_{124},v_{134},v_{234})$ with $v_{ijk}\in P_{I_i\cap I_j\cap I_k},1\leq i<j<k\leq 4$ and satisfies
\begin{equation}v_{123}-v_{124}+v_{134}-v_{234}=0.\end{equation}

The image of $\delta_1$ is generated by all the elements of the forms $\delta_1(v_{ij})$ with $v_{ij}\in P_{I_i\cap I_j},1\leq i,j\leq 4$.

\begin{center}
$\left(
  \begin{array}{cccc}
    123 & 124 & 134 & 234  \\
    v_{12} & v_{12} & 0 & 0   \\
    -v_{13} & 0 & v_{13} & 0  \\
    0 & -v_{14} & -v_{14} & 0   \\
    v_{23} & 0 & 0 & v_{23}  \\
    0 & v_{24} & 0 & -v_{24} \\
    0 & 0 & v_{34} & v_{34} \\
  \end{array}
\right)$
\end{center}

We compute $H^2(E_1)$ by constructing an epimorphism $$\psi_1: H^2(E_1)\to P_{I_1\cap I_2\cap I_3}\cap (P_{I_1\cap I_2\cap I_4}+P_{I_1\cap I_3\cap I_4}+P_{I_2\cap I_3\cap I_4})-(P_{I_1\cap I_2}+P_{I_1\cap I_3}+P_{I_2\cap I_3})$$
$$ [(v_{123},v_{124},v_{134},v_{234})]\to v_{123}+P_{I_1\cap I_2}+P_{I_1\cap I_3}+P_{I_2\cap I_3}$$
Then as a vector space $H^2(E_1)\cong \mathrm{im}\psi_1\oplus \ker \psi_1$.

For $[(v_{123},v_{124},v_{134},v_{234})]\in \ker\psi_1$, there are elements $v_{12}\in P_{I_1\cap I_2}$, $v_{13}\in P_{I_1\cap I_3}$ and $v_{23}\in P_{I_2\cap I_3}$ such that $v_{123}=v_{12}-v_{13}+v_{23}$. Combining with Equation 5, we get $(v_{124}-v_{12})-(v_{134}-v_{13})+(v_{234}-v_{23})=0$. This implies that $v_{124}-v_{12}\in P_{I_1\cap I_2\cap I_4}\cap (P_{I_1\cap I_3\cap I_4}+P_{I_2\cap I_3\cap I_4})$.

We compute $\ker \psi_1$ by constructing a homomorphism $$\psi_2: \ker \psi_1\to P_{I_1\cap I_2\cap I_4}\cap (P_{I_1\cap I_3\cap I_4}+P_{I_2\cap I_3\cap I_4})-(P_{I_1\cap I_2}\cap (P_{I_1\cap I_3}+P_{I_2\cap I_3})+P_{I_1\cap I_4}+P_{I_2\cap I_4})$$
$$[(v_{123},v_{124},v_{134},v_{234})]\to v_{124}-v_{12}+P_{I_1\cap I_2}\cap (P_{I_1\cap I_3}+P_{I_2\cap I_3})+P_{I_1\cap I_4}+P_{I_2\cap I_4})$$

It is easy to see that the homomorphism $\psi_2$ is well defined. Note that we add the item $P_{I_1\cap I_2}\cap (P_{I_1\cap I_3}+P_{I_2\cap I_3})$ to eliminate the effect of indeterminacy of $v_{12}$. By definition, $\psi_2$ is an epimorphism. So as a vector space $\ker \psi_1\cong \mathrm{im}\psi_2\oplus \ker \psi_2$.

For $[(v_{123},v_{124},v_{134},v_{234})]\in \ker\psi_2$, a similar discussion as above shows that there are elements $v_{12}\in P_{I_1\cap I_2}$, $v_{13}\in P_{I_1\cap I_3}$ and $v_{23}\in P_{I_2\cap I_3}$ with $v_{123}=v_{12}-v_{13}+v_{23}$. And there are $\gamma\in P_{I_1\cap I_2}\cap (P_{I_1\cap I_3}+P_{I_2\cap I_3})$, $v_{14}\in P_{I_1\cap I_4}$ and $v_{24}\in P_{I_2\cap I_4}$ such that $v_{124}-v_{12}=\gamma-v_{14}+v_{24}$. Similarly we can assume $\gamma=0$ if we choose $v_{12}$ and $v_{13}$ carefully in the previous step($v'_{12}=v_{12}+\gamma$ and $v'_{13}=v_{13}+\gamma$). Combining with Equation 5, we get $v_{134}-v_{13}+v_{14}=v_{234}-v_{23}+v_{24}$. This shows that $v_{134}-v_{13}+v_{14}\in P_{I_1\cap I_3\cap I_4}\cap P_{I_2\cap I_3\cap I_4}=P_{(I_1\cup I_2)\cap I_3\cap I_4}$.

Now we can construct an epimorphism $$\psi_3: \ker\psi_2\to P_{(I_1\cup I_2)\cap I_3 \cap I_4}-(P_{(I_1\cup I_2)\cap I_3}+P_{(I_1\cup I_2)\cap I_4}+P_{I_3\cap I_4})$$
$$[(v_{123},v_{124},v_{134},v_{234})]\mapsto v_{134}-v_{13}+v_{14}+P_{(I_1\cup I_2)\cap I_3}+P_{(I_1\cup I_2)\cap I_4}+P_{I_3\cap I_4}$$

We add the items $P_{(I_1\cup I_2)\cap I_3}$ and $P_{(I_1\cup I_2)\cap I_4}$ to eliminate the effect of indeterminacy of $v_{13}$ and $v_{14}$.
For $[(v_{123},v_{124},v_{134},v_{234})]\in \ker \psi_3$ we have $v_{123}=v_{12}-v_{13}+v_{23}$, $v_{124}=v_{12}-v_{14}+v_{24}$ and $v_{134}=v_{13}-v_{14}+\epsilon +\theta+v_{34}$. If we choose $v'_{13}=v_{13}+\epsilon, v'_{14}=v_{14}+\epsilon$ and $v'_{23}=v_{23}+\theta, v'_{24}=v_{24}+\theta$ in the previous step, then we have $v_{134}=v_{13}-v_{14}+v_{34}$ and $v_{234}=v_{23}-v_{24}+v_{34}$. This shows that $\ker \psi_3=0$.

Similar to the result for $H^1(E_1)$, we get

$H^2(E_1)\cong \mathrm{im}\psi_1\oplus \mathrm{im}\psi_2\oplus \mathrm{im} \psi_3$

$\cong P_{(I_1\cup I_2)\cap I_3 \cap I_4}-(P_{(I_1\cup I_2)\cap I_3}+P_{(I_1\cup I_2)\cap I_4}+P_{I_3\cap I_4})$

$\oplus P_{I_1\cap I_2\cap I_4}\cap (P_{I_1\cap I_3\cap I_4}+P_{I_2\cap I_3\cap I_4})-(P_{I_1\cap I_2}\cap (P_{I_1\cap I_3}+P_{I_2\cap I_3})+P_{I_1\cap I_4}+P_{I_2\cap I_4})$

$\oplus P_{I_1\cap I_2\cap I_3}\cap (P_{I_1\cap I_2\cap I_4}+P_{I_1\cap I_3\cap I_4}+P_{I_2\cap I_3\cap I_4})-(P_{I_1\cap I_2}+P_{I_1\cap I_3}+P_{I_2\cap I_3})$

{\bf Step 4:} For $H^3(E_1)$, it is easy to see that

$H^3(E_1)\cong P_{I_1\cap I_2\cap I_3\cap I_4}-(P_{I_1\cap I_2\cap I_3}+P_{I_1\cap I_2\cap I_4}+P_{I_1\cap I_3\cap I_4}+P_{I_2\cap I_3\cap I_4})$.

As a result, the rational cohomology group $H^*(BG(A))$ is isomorphic to $\Sigma^3 H^3(E_1)\oplus \Sigma^2 H^2(E_1) \oplus \Sigma H^1(E_1) \oplus H^0(E_1)$.
\end{example}

\begin{remark}
We will get different expressions of $H^*(BG(A))$ if we change the order of indices $I_1,I_2,\cdots,I_d$.
\end{remark}

The computation for the general case $X=X_{I_1} \cup X_{I_2}\cup \cdots\cup X_{I_d}$ is very complicate. We have the following inductive algorithm to write down the final results.

{\bf Step 1:} We begin with the cohomology group of $X_1=X_{I_1}$ which is isomorphic to $P_{I_1}$ as a graded abelian group.

{\bf Step 2:} We compute the cohomology group of $X_2=X_{I_1}\cup X_{I_2}$. The intersection $X_{I_1}\cap X_{I_2}$ has the cohomology groups $P_{I_1\cap I_2}$.
We need the following form of the M-V exact sequence which was used in Zhao and Gao\cite{Zhao_Gao_20}
\begin{lemma}
Let $X=X_1\cup_{X_0} X_2$ be the push-out of the diagram $X_1 \stackrel{j_1} \leftarrow X_0\stackrel{j_2}\rightarrow X_2$. The homomorphism $j: H^*(X_1)\oplus H^*(X_2)\to H^*(X_0)$ is given by $j(u,v)=j_1^*(u)-j_2^*(v)$. If $X_1,X_2$ are deformation retracts of some open subspaces of $X$, then there exists a short exact sequence \begin{center}$0\to \Sigma \mathrm{coker}j\to H^*(X)\to \ker j\to 0.$\end{center} \label{lm2.6}
\end{lemma}
Applying this lemma to $X_{I_1}\cup X_{I_2}$, there are two homomorphisms $j_1: H^*(X_{I_1})\rightarrow H^*(X_{I_1\cap I_2}), $ and $j_2: H^*(X_{I_2})\rightarrow H^*(X_{I_1\cap I_2})$. Then $\ker j\cong P_{I_1 \cup I_2}$ and $\mathrm{coker} j\cong P_{I_1\cap I_2}-(P_{I_1}+P_{I_2})$. And we get
\begin{center}$H^*(X_1\cup X_2)\cong \Sigma \mathrm{coker} j\oplus \ker j\cong \Sigma(P_{I_1\cap I_2}-(P_{I_1}+P_{I_2}))\oplus P_{I_1 \cup I_2}$.\end{center}
We denote the two items in the expression by $F_1$ and $F_2$. Then $F_1=\Sigma(P_{I_1\cap I_2}-(P_{I_1}+P_{I_2}))$ has the form $\Sigma(A-B)$ and $F_2$ has the form $A'$.

{\bf Step 3:} We compute the cohomology group of $X_3=X_{I_1}\cup X_{I_2}\cup X_{I_3}$. Note that $X_3=X_2 \cup X_{I_3}$ and the intersection of $X_2$ and $X_{I_3}$ is $\bar X_2=(X_{I_1}\cup X_{I_2})\cap X_{I_3}=X_{I_1\cap I_3}\cup X_{I_2\cap I_3}$. By applying the result of Step 2 to $\bar I_1=I_1\cap I_3$ and $\bar I_2= I_2\cap I_3$, we get

$H^*(X_{I_1\cap I_3} \cup X_{I_2\cap I_3})\cong \Sigma \mathrm{coker} \bar j\oplus \ker \bar j \cong \Sigma (P_{(I_1\cap I_3)\cap (I_2\cap I_3)}-(P_{I_1\cap I_3}+P_{I_2\cap I_3}))\oplus P_{(I_1\cap I_3)\cup (I_2\cap I_3)}$.

And we have two homomorphisms
$j_1: H^*(X_{I_1}\cup X_{I_2})\rightarrow H^*(X_{(I_1\cap I_3)\cup (I_2\cap I_3)})$ and $j_2: H^*(X_{I_3})\rightarrow H^*(X_{(I_1\cap I_3)\cup (I_2\cap I_3)})$.
These homomorphisms can be depicted in the following diagram.

\noindent $\xymatrix{
  H^*(X_{I_1}\cup X_{I_2}) \ar[d]^{j_1} &\cong & \Sigma(P_{I_1\cap I_2}-(P_{I_1}+P_{I_2}))  \ar[d]^{j_1^{\mathrm{coker}}}& \oplus & P_{I_1 \cup I_2}  \ar[d]^{j_1^{\ker }}\\
  H^*(X_{(I_1\cap I_3)\cup (I_2\cap I_3)})&\cong & \Sigma P_{(I_1\cap I_3)\cap (I_2\cap I_3)}-(P_{I_1\cap I_3}+P_{I_2\cap I_3}) &\oplus & P_{(I_1\cap I_3)\cup (I_2\cap I_3)} \\
  H^*(X_{I_3}) \ar[u]^{j_2}&\cong& 0 \ar[u]^{j_2^{\mathrm{coker}}}&\oplus & P_{I_3}  \ar[u]^{j_2^{\ker }}}$

In the M-V exact sequence we combine the homomorphisms $j_1,j_2$ to give the homomorphism $j$ and $H^*(X_{I_1}\cup  X_{I_2}\cup X_{I_3})\cong \Sigma \mathrm{coker} j\oplus \ker j$. If we assume that the homomorphisms $j_1$ and $j_2$ are decomposed into sum of its summand, then from the M-V exact sequence we get $H^*(X_{I_1}\cup  X_{I_2}\cup X_{I_3})\cong \Sigma \mathrm{coker} j_1^{\mathrm{coker}} \oplus \Sigma \mathrm{coker} j_2^{\ker} \oplus \ker j_1^{\mathrm{coker}} \oplus \ker j_2^{\ker}$. Detailed computation gives the following rule to obtain $H^*(X_3)$ from $H^*(X_2)$.
%The $4$-item as in the above algorithm originate this way.

%By direct computation we will get the same result as described before.

%Step4 is the same but more complex.

% whose cohomology can be computed similarly by Step 2. Now We compute the cohomology of $X_3$ still by using M-V exact sequence. We describe the algorithm to compute $H^*(X_3)$.

i) We write the cohomology of $\bar X_2=X_{I_1\cap I_3}\cup X_{I_2\cap I_3}$ below the cohomology of $X_2=X_{I_1} \cup X_{I_2}$ and we get

$\Sigma(P_{I_1\cap I_2}-(P_{I_1}+P_{I_2}))\ \ \ \ \ \ \ \ \ \ \ \oplus P_{I_1\cup I_2}$

$\Sigma(P_{I_1\cap I_2\cap I_3}-(P_{I_1\cap I_3}+P_{I_2\cap I_3}))\oplus P_{(I_1\cup I_2)\cap I_3}$

The result for $\bar X_2$ is obtained from the result for $X_2$ by replacing $I_1$ and $I_2$ to $\bar I_1=I_1\cap I_3$ and $\bar I_2=I_2\cap I_3$.
The two expressions contain four direct sum items, that is $\Sigma(P_{I_1\cap I_2}-(P_{I_1}+P_{I_2}))$, $P_{I_1\cup I_2}$, $\Sigma(P_{I_1\cap I_2\cap I_3}-(P_{I_1\cap I_3}+P_{I_2\cap I_3}))$ and $P_{(I_1\cup I_2)\cap I_3}$.

ii) Each of these items will be transformed into an item in the expression of $H^*(X_3)$.
The two corresponding items $\Sigma(P_{I_1\cap I_2}-(P_{I_1}+P_{I_2}))$ and $\Sigma(P_{I_1\cap I_2\cap I_3}-(P_{I_1\cap I_3}+P_{I_2\cap I_3}))$ have the forms as $\Sigma(A-B)$ and $\Sigma(\bar A-\bar B)$ with $B\subset A,\bar B\subset \bar A$ and $A\subset \bar A, B\subset \bar B$. They will be transformed into two new items $\Sigma(A\cap \bar B-B) $ and $\Sigma^2(\bar A-(\bar B+A)$. The two items $P_{I_1\cup I_2}$ and $P_{(I_1\cup I_2)\cap I_3}$ have the forms $A'$ and $\bar A'$ with $A'\subset \bar A'$, which are transformed into the two items $\Sigma(\bar A'-(A'+P_{I_3}))$ and $A'\cap P_{ I_3}$. We write the four items as $F_1\oplus F_2\oplus F_3\oplus F_4$, and require $F_4$ to be the unique item without grade lifting(no $\Sigma$ appears).

And we continue this process until {\bf Step $k$}.

Suppose we have finished the step $k$ and obtained the cohomology group $H^*(X_{k})=H^*(X_{I_1}\cup X_{I_2}\cup \cdots \cup X_{I_{k}})$ and it has the form
of $F_1\oplus F_2 \oplus \cdots \oplus F_{2^{k-1}-1}\oplus F_{2^{k-1}}$. And we require that there is no $\Sigma$ in $F_{2^{k-1}}$.

{\bf Step $k+1$:} now we write $H^*(\bar X_{k})=H^*(X_{I_1\cap I_{k+1}}\cup X_{I_2\cap I_{k+1}}\cup \cdots \cup X_{I_{k}\cap I_{k+1}})$ below the cohomology groups $H^*(X_{k})$ and we get

$F_1\oplus F_2 \oplus \cdots \oplus F_{2^{k-1}-1}\oplus F_{2^{k-1}}$

$\bar F_1\oplus \bar F_2 \oplus \cdots \oplus \bar F_{2^{k-1}-1}\oplus  \bar F_{2^{k-1}}$

Through the procedure similar to that in step 3 we know that the two items $F_i$ and $\bar F_i, 1\leq i\leq 2^{k-1}-1$ have the forms of $\Sigma^r(A-B)$ and $\Sigma^r(\bar A-\bar B)$ with $B\subset A,\bar B\subset \bar A$ and $A\subset \bar A, B\subset \bar B$, which will be transformed into two new items $\Sigma^r(A\cap \bar B-B)$ and $\Sigma^{r+1}(\bar A-(\bar B+A))$ in $H^*(X_{k+1})$. The last two items $F_{2^{k-1}}$ and $\bar F_{2^{k-1}}$ have the forms $A'$ and $\bar A'$ with $A'\subset \bar A'$, which are transformed into the two items $\Sigma(\bar A'-(A'+P_{I_{k+1}}))$ and $A' \cap P_{I_{k+1}}$.

%The algorithm is highlighted by the following inaccurate computation procedure.

The algorithm ends when $k=d-1$.\\

The following is an example illustrating the algorithm for $d=4$.
\begin{example}
{\bf Step 1}. $X_{I_1}$, input $P_{I_1}$

{\bf Step 2}. $X_{I_1}\cup X_{I_2}$

We write the cohomology of $X_1$ and $\bar X_1$.

$P_1$

$P_{I_1\cap I_2}$

And we get $H^*(X_2)\cong \Sigma(P_{I_1\cap I_2}-(P_{I_1}+P_{I_2}))\oplus P_{I_1\cup I_2}$

{\bf Step 3}. $X_{I_1}\cup X_{I_2}\cup X_{I_3}$

$(X_{I_1}\cup X_{I_2})\cap X_{I_3}=(X_{I_1}\cap X_{I_3}) \cup (X_{I_2}\cap X_{I_3})=X_{I_1\cap I_3}\cup X_{I_2\cap I_3}$

We write the cohomology of $X_2$ and $\bar X_2$.

$\Sigma(P_{I_1\cap I_2}-(P_{I_1}+P_{I_2}))\oplus P_{I_1\cup I_2}$

$\Sigma(P_{I_1\cap I_2\cap I_3}-(P_{I_1\cap I_3}+P_{I_2\cap I_3}))\oplus P_{(I_1\cup I_2)\cap I_3}$

And we get $H^*(X_3)\cong \Sigma^2(P_{I_1\cap I_2\cap I_3}-(P_{I_1\cap I_2}+P_{I_1\cap I_3}+P_{I_2\cap I_3}))\oplus \Sigma(P_{I_1\cap I_2}\cap (P_{I_1\cap I_3}+P_{I_2\cap I_3})-(P_{I_1}+P_{I_2}))\oplus \Sigma(P_{(I_1\cup I_2)\cap I_3}-(P_{I_1\cup I_2}+P_{I_3})) \oplus P_{I_1\cup I_2\cup I_3}$

{\bf Step 4}. $X_{I_1}\cup X_{I_2}\cup X_{I_3}\cup X_{I_4}$

$(X_{I_1}\cup X_{I_2}\cup X_{I_3})\cap X_{I_4}=X_{I_1\cap I_4} \cup X_{I_2\cap I_4}\cup X_{I_3\cap I_4}$

We write the cohomology of $X_3$ and $\bar X_3$.

$\Sigma^2(P_{I_1\cap I_2\cap I_3}-(P_{I_1\cap I_2}+P_{I_1\cap I_3}+P_{I_2\cap I_3}))\oplus \Sigma(P_{I_1\cap I_2}\cap (P_{I_1\cap I_3}+P_{I_2\cap I_3})-(P_{I_1}+P_{I_2}))\oplus \Sigma(P_{(I_1\cup I_2)\cap I_3}-(P_{I_1\cup I_2}+P_{I_3})) \oplus P_{I_1\cup I_2\cup I_3}$

$\Sigma^2(P_{I_1\cap I_2\cap I_3\cap I_4}-(P_{I_1\cap I_2\cap I_4}+P_{I_1\cap I_3\cap I_4}+P_{I_2\cap I_3\cap I_4}))\oplus \Sigma(P_{I_1\cap I_2\cap I_4}\cap (P_{I_1\cap I_3\cap I_4}+P_{I_2\cap I_3\cap I_4})-(P_{I_1\cap I_4}+P_{I_2\cap I_4}))\oplus \Sigma(P_{(I_1\cup I_2)\cap I_3 \cap I_4}-(P_{(I_1\cup I_2)\cap I_4}+P_{I_3\cap I_4})) \oplus P_{(I_1\cup I_2\cup I_3)\cap I_4}$

And we get

$H^*(X_4)\cong \Sigma^3(P_{I_1\cap I_2\cap I_3\cap I_4}-(P_{I_1\cap I_2\cap I_3}+P_{I_1\cap I_2\cap I_4}+P_{I_1\cap I_3\cap I_4}+P_{I_2\cap I_3\cap I_4}))$

$\oplus \Sigma^2(P_{(I_1\cup I_2)\cap I_3 \cap I_4}-(P_{(I_1\cup I_2)\cap I_3}+P_{(I_1\cup I_2)\cap I_4}+P_{I_3\cap I_4}))$

$\oplus \Sigma^2(P_{I_1\cap I_2\cap I_4}\cap (P_{I_1\cap I_3\cap I_4}+P_{I_2\cap I_3\cap I_4})-(P_{I_1\cap I_2}\cap (P_{I_1\cap I_3}+P_{I_2\cap I_3})+P_{I_1\cap I_4}+P_{I_2\cap I_4}))$

$\oplus \Sigma^2(P_{I_1\cap I_2\cap I_3}\cap (P_{I_1\cap I_2\cap I_4}+P_{I_1\cap I_3\cap I_4}+P_{I_2\cap I_3\cap I_4})-(P_{I_1\cap I_2}+P_{I_1\cap I_3}+P_{I_2\cap I_3}))$

$\oplus \Sigma(P_{(I_1\cup I_2\cup I_3)\cap I_4}-(P_{I_1\cup I_2\cup I_3}+P_{I_4}) $

$\oplus \Sigma(P_{(I_1\cup I_2)\cap I_3}\cap (P_{(I_1\cup I_2)\cap I_4}+P_{I_3\cap I_4})-(P_{I_1\cup I_2}+P_{I_3}))$

$\oplus \Sigma(P_{I_1\cap I_2}\cap (P_{I_1\cap I_3}+P_{I_2\cap I_3})\cap (P_{I_1\cap I_4}+P_{I_2\cap I_4})-(P_{I_1}+P_{I_2}))$

$\oplus P_{I_1\cup I_2\cup I_3\cup I_4}$
\end{example}

By classifying the maximal objects in $C(A)$ and using the above algorithm, we obtain the results for rational cohomology groups of the classifying spaces of Kac-Moody groups of infinite type with $n=3$ and $4$.

\begin{example}
The rational cohomology groups of the classifying space $BG(A)$ for $n=3$ are listed case by case as below.

1. $\{1\},\{2\},\{3\}$

$\Sigma(P-(P_{1}+P_{2}))\oplus \Sigma(P-(P_{12}+P_{3}))\oplus P_{123}$

2. $\{1,2\}, \{3\}$ \ \  \ \

$\Sigma(P-(P_{12}+P_{3}))\oplus P_{123}$

3. $\{1,2\}, \{1,3\}$  \ \ \ \

$\Sigma(P_1-(P_{12}+P_{13}))\oplus P_{123}$

4. $\{1,2\},\{1,3\}, \{2,3\}$

$\Sigma^2(P-(P_{1}+P_{2}+P_3))\oplus \Sigma(P_1\cap(P_2+P_3)-(P_{12}+P_{13}))\oplus P_{123}$
\end{example}

\begin{example}
The rational cohomology groups of the classifying space $BG(A)$ for $n=4$ are listed case by case as below.

1. $\{1\},\{2\},\{3\},\{4\}$, \ \  \ \ \ \

$\Sigma(P-(P_{1}+P_{2}))\oplus \Sigma(P-(P_{12}+P_{3}))\oplus \Sigma(P-(P_{123}+P_{4}))\oplus P_{1234}$

2. $\{1,2\}, \{3\},\{4\}$, \ \  \ \ \ \

$\Sigma(P-(P_{12}+P_{3}))\oplus \Sigma(P-(P_{123}+P_{4}))\oplus P_{1234}$

3. $\{1,2\}, \{3,4\}$,

 $\Sigma(P-(P_{12}+P_{34}))\oplus P_{1234}$

4. $\{1,2\}, \{1, 3\}, \{4\}$

$\Sigma(P_1-(P_{12}+P_{13}))\oplus \Sigma(P-(P_{123}+P_{4}))\oplus P_{1234}$

5. $\{1,2\}, \{1, 3\}, \{1,4\}$

$\Sigma(P_1-(P_{12}+P_{13}))\oplus \Sigma(P_1-(P_{123}+P_{14}))\oplus P_{1234}$

6. $\{1,2\}, \{1,3\},\{3,4\},$

$\Sigma(P_1-(P_{12}+P_{13}))\oplus \Sigma(P_3-(P_{123}+P_{34}))\oplus P_{1234}$

7. $\{1,2\}, \{1,3\}, \{2,3\}, \{4\}$

$\Sigma^2(P-(P_1+P_{2}+P_{3}))\oplus \Sigma (P_1\cap (P_2+P_3)-(P_{12}+P_{13})) \oplus \Sigma(P-(P_{123}+P_{4}))\oplus P_{1234}$

8. $\{1,2\}, \{1,3\},\{1,4\},\{2,3\}$

$\Sigma^2(P-(P_1+P_2+P_3))\oplus \Sigma (P_1\cap (P_2+P_3)-(P_{12}+P_{13}))\oplus \Sigma(P_1-(P_{123}+P_{14}))\oplus P_{1234}$

9. $\{1,2\}, \{1,3\}, \{2,4\},\{3,4\}$

$\Sigma^2(P-(P_2+P_3+P_4))\oplus \Sigma(P_1-(P_{12}+P_{13}))\oplus \Sigma(P_{2}\cap (P_3+P_4)-(P_{123}+P_{24}))\oplus P_{1234}$

10. $\{1,2\},\{1,3\},\{1,4\},\{2,3\},\{2,4\}$

$\Sigma^2(P-(P_1+P_2+P_3))\oplus \Sigma^2(P-(P_{1}+P_2+P_{4}))\oplus \Sigma(P_1\cap(P_2+P_4)-(P_{123}+P_{14})\oplus \Sigma(P_1\cap(P_2+P_3)-(P_{12}+P_{13})\oplus P_{1234}$

11. $\{1,2\},\{1,3\},\{1,4\},\{2,3\},\{2,4\},\{3,4\}$

$\Sigma^2(P-(P_1+P_2+P_3))\oplus \Sigma^2(P-(P_{1}+P_2+P_{4}))\oplus \Sigma^2(P-(P_1\cap(P_2+P_4)+P_3+P_4))\oplus \Sigma(P_1\cap(P_2+P_4)\cap(P_3+P_4)-(P_{123}+P_{14})\oplus \Sigma(P_1\cap(P_2+P_3)-(P_{12}+P_{13})\oplus P_{1234}$

12. $\{1,2,3\}, \{4\}$

$\Sigma(P-(P_{123}+P_{4}))\oplus P_{1234}$

13. $\{1,2,3\}, \{1,4\}$

$\Sigma(P_1-(P_{123}+P_{14}))\oplus P_{1234}$

14. $\{1,2,3\},\{1,4\},\{2,4\},$

$\Sigma^2(P-(P_{1}+P_{2}+P_4))\oplus \Sigma(P_{1}\cap(P_2+P_4)-(P_{123}+P_{14}))\oplus P_{1234}$

15. $\{1,2,3\}, \{1,4\},\{2,4\}, \{3,4\}$

$\Sigma^2(P-(P_{1}+P_{2}+P_4))\oplus \Sigma^2(P-(P_{1}\cap (P_2+P_4)+P_{3}+P_4))\oplus \Sigma((P_{1}\cap (P_2+P_4)\cap (P_3+P_4)-(P_{123}+P_{14})\oplus P_{1234}$

16. $\{1,2,3\}, \{1,2,4\}$

$\Sigma(P_{12}-(P_{123}+P_{124}))\oplus P_{1234}$

17. $\{1,2,3\}, \{1,2,4\}, \{3,4\}$

$\Sigma^2(P-(P_{12}+P_3+P_4))\oplus \Sigma(P_{12}\cap (P_3+P_4)-(P_{123}+P_{124})) \oplus P_{1234}$

18. $\{1,2,3\}, \{1,2,4\}, \{1,3,4\}$

$\Sigma^2(P_{1}-(P_{12}+P_{13}+P_{14}))\oplus \Sigma(P_{12}\cap (P_{13}+P_{14})-(P_{123}+P_{124}))\oplus P_{1234}$.

19. $\{1,2,3\}, \{1,2,4\},\{1,3,4\}, \{2,3,4\}$

$\Sigma^3(P-(P_{1}+P_{2}+P_{3}+P_{4}))\oplus \Sigma^2(P_{1}\cap (P_{2}+P_{3}+P_{4})-(P_{12}+P_{13}+P_{14}))\oplus \Sigma^2(P_{2}\cap (P_{3}+P_{4})-(P_{12}\cap (P_{13}+P_{14})+P_{23}+P_{24}))\oplus \Sigma(P_{12}\cap (P_{13}+P_{14})\cap (P_{23}+P_{24})-(P_{123}+P_{124}))\oplus P_{1234}$.

%11,15,17,18,19 multiplication?
\end{example}

By the above algorithm there are always $2^{d-1}$ items in the cohomology of $H^*(BG(A))$. But some items are zero for trivial reason. We delete them from the lists. We have done the similar computations for $n=5$ with the help of Zhang Yukun, Luo Yajing, Gong Tao, Wang Ran, Li Yingxin, Qu Qingrui, Liang Peng and Gu Peiqi. The result is too lengthy to be listed here.

%\begin{example}14,15,19 multiplication structure.\end{example}

%\begin{define}number of simplexes\end{define}

\section{The polynomial invariants of Weyl groups and their subgroups}
To compute the rational cohomology of the classifying space $BG(A)$, it is vital to know the ring of polynomial invariants $P_I$ for $I\subset S$(see the expressions in Examples 4.3 and 4.4).  In some cases the computation of $P_I$ for $I\subset S$ can be reduced to the computation of rational polynomial invariants of Weyl group $W(A_I)$ for the Cartan matrix $A_I$.

%For example, the affine case is determined by Wang\cite{Wang_19} and the indefinite case by Zhao and Jin\cite{Zhao_Jin_14}.

Each decomposable $n\times n$ Cartan matrix $A$ can be written, up to order of the elements in $S$, as a sum $A=\bigoplus\limits_{i=1}^k A_{S_i}$ of sub-matrices with $S$ the disjoint union of $S_i$. And the corresponding Kac-Moody groups can be decomposed into
$G(A)=\prod \limits_{i=1}^k G(A_{S_i})$. For each $I\subset S$ the parabolic subgroup $G_I(A)$ has a decomposition $G_I(A)=\prod\limits_{i=1}^k G_{I\cap S_i}(A_i)$. And the Weyl group of $G_I(A)$ is isomorphic to $\prod\limits_{i=1}^k W_{I\cap S_i} (A_i)$. It is obvious that the ring of invariants $P_I$ is isomorphic to $\bigotimes\limits_{i=1}^k P_{I\cap S_i}$. So the computation of $P_I$ can be reduced to the indecomposable case.

%\begin{lemma}For an indecomposable $n\times n$ Cartan matrix $A$ and $I\subset S$, then $P_I\cong $ the tensor product of the $W(A_I)$ rational polynomial invariants ring of $\Q[w_1,w_2,\cdots,w_n]$ and $\Q[w_{|I|+1},w_{|I|+2},\cdots, w_n]$\end{lemma}

Let $A=(a_{ij})_{n\times n}$ be a Cartan matrix. The Weyl reflections $\sigma_1,\sigma_2,\cdots,\sigma_n$ satisfy
$$\sigma_i(w_j)=\left\{
                                                     \begin{array}{ll}
                                                       w_j-\alpha_j, & \hbox{ for } i=j; \\
                                                       w_j, & \hbox{ for } i\not=j.
                                                     \end{array}
                                                   \right.$$
Thus we have $$(\sigma_1(w_1),\sigma_2(w_2),\cdots,\sigma_n(w_n))=(w_1,w_2,\cdots,w_n)(I_n-A).$$

We write $(w_1,w_2,\cdots,w_n)$ as $w$. For $S=I\sqcup J$, we write $w$ as $(w_I,w_J)$, where $w_I=(w_1,w_2,\cdots,w_{|I|})$ and $w_J=(w_{|I|+1},w_{|I|+2},\cdots,w_n)$. We have $$(\sigma_1(w_1),\sigma_2(w_2),\cdots,\sigma_{|I|}(w_{|I|}))=(w_1,w_2,\cdots,w_n)\left(
                                                                          \begin{array}{c}
                                                                            I_{|I|}-A_{I} \\
                                                                            -A_{JI} \\
                                                                          \end{array}
                                                                        \right)$$
and $\sigma_i(\omega_j)=\omega_j$ for all $i\in I, j\in J$, where $A_{JI}$ is the submatrix $A_{JI}=(a_{ji})_{j \in J, i\in I}$ of $A$.

By choosing a $|J|\times |I|$ matrix $C$ with rational elements, we get another basis for the rational linear space of weights, $$w'=(w'_I,w'_J)=(w_I,w_J)
\left(\begin{array}{cc}
                                                                            I_{|I|} & 0 \\
                                                                            C & I_{|J|} \\
                                                                          \end{array}
                                                                        \right). $$
By a direct computation we get
$$(\sigma_1(w'_1),\sigma_2(w'_2),\cdots,\sigma_{|I|}(w'_{|I|}))=(w_1,w_2,\cdots,w_n)\left(
                                                                                   \begin{array}{c}
                                                                                     I_{|I|} -A_{I}\\
                                                                                     C-A_{JI} \\
                                                                                   \end{array}
                                                                                 \right)
=(w'_1,w'_2,\cdots,w'_n)\left(
                                                                                   \begin{array}{c}
                                                                                     I_{|I|} -A_{I}\\
                                                                                     C A_I -A_{JI} \\
                                                                                   \end{array}
                                                                                 \right)$$
If the rank of matrix $A_I$ is the same as that of $\left(
                                                                                   \begin{array}{c}
                                                                                     A_{I}\\
                                                                                     A_{JI} \\
                                                                                   \end{array}
                                                                                 \right)$,
then the rational matrix $C$ can be chosen such that $CA_I=A_{JI}$. As a consequence we get $$(\sigma_1(w'_1),\sigma_2(w'_2),\cdots,\sigma_{|I|}(w'_{|I|}))=(w'_1,w'_2,\cdots,w'_{|I|})
                                                                                   (\begin{array}{c}
                                                                                     I_{|I|} -A_{I}
                                                                                   \end{array}).$$
This implies that the linear action of $W_I(A)$ on the $\Q$-linear space spanned by $w'_1,w'_2,\cdots,w'_{|I|}$ is isomorphic to the action of $W(A_I)$ on the $\Q$-linear space of weights for $G(A_I)$. Combining with the fact that $W_I(A)$ acts trivially on $w_j, j\in J$, we prove the following proposition.

\begin{prop}
Let $A$ be an $n\times n$ Cartan matrix and the matrix equation $XA_I=A_{JI}$ has a rational solution $C$, then the ring of invariants $P_I$ of $W_I(A)$ is isomorphic to the tensor product of the ring of invariants of the Weyl group $W(A_I)$ acting on $\Q$-linear space of weights for $G(A_I)$ and the polynomial ring $\Q[w_{|I|+1}, w_{|I|+2}, \cdots, w_{n}]$.
\end{prop}

Particularly the proposition is true if $\det A_I\not=0$.

Even one knows all the rings of invariants $P_I$, it is still difficult to compute the rational cohomology groups of $BG(A)$ explicitly. There are expressions like $P-(P_1+P_2+P_3)$ and $P_{1}\cap (P_2+P_4)\cap (P_3+P_4)-(P_{123}+P_{14})$ for which we even can not compute their Poincar\'{e} series.

In \cite{Foley_15} Foley claimed that for a $3\times 3$ Cartan matrix $A$ of infinite type, $P-(P_1+P_2+P_3)=0$. This result hints the following conjecture.

\begin{conj}
For an $n\times n$ Cartan matrix $A$ of infinite type, $P=P_1+P_2+\cdots+P_n$.
\end{conj}

%We further conjecture that for a Cartan matrix of indefinite tyoe all the items in Example 4.4 contains $\Sigma^2$ or $\Sigma^3$ are trivial. And this result is unlikely true for $n\geq 5$ generally.

\section{Some applications and conjectures}
%At first we note that for any coefficient ring the M-V spectral sequence converges to $H^*(BG(A))$. But in general it does not collapse.

The computation for the cohomology groups of $BG(A)$(or $BK(A)$) has several applications.

\begin{theorem}
For a Kac-Moody group $K(A)$ of infinite type, the image of the Borel homomorphism $Bi^*: H^*(BK(A);\Q)\rightarrow H^*(BT,\Q)$ is the ring of invariants $P^{W(A)}$.
\end{theorem}

\noindent {\bf Proof: }We have $BK(A)\sim X_{I_1}\cup X_{I_2} \cup \cdots \cup X_{I_d}$. We can also write $BT$ as $BT\cup BT\cup \cdots \cup BT$ and there is a
homomorphism $f$ between the M-V spectral sequences of $H^*(BK(A))$ and $H^*(BT)$. Borel's theorem asserts that the homomorphism $H^*(X_{I_i})\rightarrow H^*(BT)$ has image the ring of invariants $P^{W_{I_i}(A)}$. Since $S=\bigcup\limits_{i=1}^d {I_i}$, the image of $Bi^*$ is isomorphic to $\bigcap\limits_{i=1}^d P^{W_{I_i}(A)}=P^{W(A)}$.

%\begin{theorem}
%For a Kac-Moody group of infinite type the image of the Borel homomorphism $Bi^*: H^*(K(A))\rightarrow H^*(BT)$ is the ring of invariants $P^{W(A)}$.
%\end{theorem}

%\begin{theorem}
%The mod $p$ cohomology ring of a Kac-Moody group $G(A)$ is finite generated.
%\end{theorem}

The following result is proved in Ruan and Zhao\cite{Ruan_Zhao_20}.

\begin{theorem}
For a rank $3$ Kac-Moody group $K(A)$ of infinite type, the integral cohomology $H^*(BK(A);\Z)$ has $p$-torsion elements for any prime number $p$.
\end{theorem}

Based on theorem 6.2 we can give the following conjecture.

\begin{conj}
For a Kac-Moody group $K(A)$ of infinite type, $H^*(BG(A);\Z)$ has $p$-torsion elements for any prime $p$.
\end{conj}

In \cite{Zhao_Gao_20} we show that for a generic Kac-Moody group, the rational cohomology ring $H^*(BG(A);\Q)$ is infinitely generated. So it is natural to give the following conjecture.

\begin{conj}
For a Cartan matrix $A=(a_{ij})_{n\times n}$ of indefinite type with $n\geq 3$, the rational cohomology ring $H^*(BG(A);\Q)$ is infinitely generated.
\end{conj}

\begin{remark}
For Conjecture 6.2 it is known that for Cartan matrix $A=(a_{ij})_{n\times n}$ of finite type, affine type or $n\leq 2$, the rational cohomology ring $H^*(BG(A);\Q)$ is finitely generated.
\end{remark}

The computation in this paper can be used to attack the above conjectures.

\end{document}